\documentclass[reqno,11pt]{amsart}
\usepackage{amsmath} 
\usepackage{amsthm, amsfonts, amsgen, stmaryrd,enumerate}
\usepackage[english]{babel}
\usepackage{fullpage}
\usepackage[centertags]{amsmath}
\usepackage{indentfirst}
\usepackage{fancyhdr}
\usepackage[dvips]{graphicx}
\usepackage{psfrag}
\numberwithin{equation}{section}
\usepackage[latin1]{inputenc}

\usepackage{color}

\title[Hodge dualities over $\SU$]{(A class of) Hodge duality operators \\ over the quantum SU(2)}
\date{3 October  2012}
\author{Alessandro Zampini}
\address{ Hausdorff Zentrum der Universit\"at Bonn, Endenicher Allee 62, D-53115 Bonn, Germany;  \newline\indent
Max Planck Institut f\"ur  Mathematik - Bonn, 
Vivatsgasse 7, D-53111 Bonn, Germany.
\newline\indent {\it Present Affiliation:}
\newline\indent Mathematisches Institut der L.M.U., Theresienstra\ss e 39,  D-80333 M\"unchen, Germany.}
 \email{zampini@math.lmu.de}

\newtheorem{theo}{Theorem}[section]

\newtheorem{prop}[theo]{Proposition}

\newtheorem{rema}[theo]{Remark}

\newcommand{\nn}{\nonumber}

\newcommand{\dd}{{\rm d}}
\newcommand{\ca}{\mathcal{A}}

\newcommand{\ch}{\mathcal{H}}
\newcommand{\cl}{\mathcal{L}}

\newcommand{\cq}{\mathcal{Q}}

\newcommand{\cu}{\mathcal{U}}        
\newcommand{\SU}{\mathrm{SU}_q(2)}  
\newcommand{\ASU}{\ca(\mathrm{SU}_q(2))}  
\newcommand{\sq}{\mathrm{S}^2_{q}}  
\newcommand{\Asq}{\ca(\mathrm{S}^2_{q})}  
\newcommand{\su}{\cu_q(\mathfrak{su}(2))}  

\newcommand{\eps}{\varepsilon}      
\newcommand{\co}[2]{#1_{(#2)}}      
\newcommand{\hs}[2]{\left\langle #1,#2\right\rangle}  
\newcommand{\hp}[2]{\left\{ #1,#2\right\}} 

\newcommand{\lt}{{\triangleright}}    
\newcommand{\rt}{{\triangleleft}}
\newcommand{\IC}{{\mathbb C}} 
\newcommand{\IR}{{\mathbb R}} 
\newcommand{\IN}{{\mathbb N}} 
\newcommand{\IZ}{{\mathbb Z}} 
\DeclareMathOperator{\Ad}{Ad}       
\DeclareMathOperator{\id}{id}       
\DeclareMathOperator{\U}{U}       

\newcommand{\figureheight}{8cm}
\newcommand{\putfig}[2]{\begin{figure}[htp]
        \special{isoscale c:/itex/texfig/#1.wmf, \the\hsize \figureheight}
        \vspace{\figureheight}
        \caption{#2}\label{fig:#1}
        \end{figure}}
\newcommand{\pictureheight}{4cm}
\newcommand{\putpicture}[2]{\begin{figure}[htp]
        \special{isoscale c:/itex/texfig/#1.wmf, \the\hsize \pictureheight}
        \vspace{\pictureheight}
        \caption{#2}\label{fig:#1}
        \end{figure}}

\newcommand{\beqa}{\begin{eqnarray}}
\newcommand{\eeqa}{\end{eqnarray}}
\newcommand{\beq}{\begin{equation}}
\newcommand{\eeq}{\end{equation}}

\newcommand{\ma}{\mathit{a}}

\newcommand{\mT}{\mathfrak{T}}
\newcommand{\mL}{\mathfrak{L}}
\newcommand{\mg}{\mathit{g}_I}

\begin{document}

\thispagestyle{empty}
\begin{abstract}
\noindent 
On the exterior algebra over $\SU$ coming from the four dimensional bicovariant calculus \`a la Woronowicz we introduce, using sesquilinear contraction maps, a class of metrics and Hodge duality operators, and compare this formulation with the previous ones presented in the literature.

\end{abstract}

\maketitle

\tableofcontents

\section{Introduction}
This paper is a further step along the study of the geometry that the quantum $\mathrm{SU}(2)$ group acquires, when it is equipped with the $4D_{+}$ first order  bicovariant calculus \`a la Woronowicz and the corresponding (not universal, four dimensional) exterior algebra $\Omega(\SU)$. It is aimed to evolve the analysis presented in \cite{ag10}, where, using the canonical braiding of the calculus, a definition of  Hodge operators on both $\Omega(\SU)$ and $\Omega(\sq)$ (with $\sq$ the standard Podle\'s sphere) were proposed. The squares of those Hodge operators do not show  the classical degenerate spectra: they are not constant on each $k$-form in $\Omega^{k}(\SU)$, where  they turn out to be not even necessarily diagonal.

Is it possible to introduce novel Hodge operators on $\Omega(\SU)$, whose square is diagonal, and such that its spectrum  has a degeneracy 
which naturally fits with the degeneracy of the antisymmetriser operators (used to construct the exterior algebra) coming from the canonical (non trivial) braiding associated to the $4D_{+}$ first order calculus?

To answer this question, we merge the formulation described in \cite{ag10} with (again!) that  in \cite{hec99}, and introduce a notion of Hodge operators on $\Omega(\SU)$ whose square \emph{is} diagonal, with a spectrum showing the \emph{same} degeneracy of the antisymmetrisers of the calculus. Such Hodge operators and their corresponding metric structures present close similarities to those studied in \cite{hec99}, as well as interesting differences with those in \cite{ag10} and in \cite{gae1, gae2}, where a Hodge map on $\Omega(\SU)$ has been introduced via a suitable $\epsilon$ tensor.

\subsection{Overview}
It is well known what Hodge duality operators on classical group manifolds are. Consider   an $N$-dimensional  real connected  Lie group  $G$. 
A metric on the group $G$ is a non degenerate tensor $g:\mathfrak{X}(G)\otimes\mathfrak{X}(G)\to\ca(G)$  (with $\ca(G)={\rm Fun}(G)$ the $*$-algebra of complex valued coordinate functions on $G$)  which is symmetric, that is $g(X,Y)=g(Y,X)$ for any pair of vector fields $X,Y\in\,\mathfrak{X}(G)$,  and real, $g^*(X,Y)=g(Y^*,X^*)$. Written in its normal form on  a real basis $\{\theta^{a}=\theta^{a*}, a=1, \dots, N\}$ of $\Omega^1(G)$, it is
\beq
\label{gmf}
g=\sum_{a,b=1}^{N}\eta_{ab}\,\theta^{a}\otimes\theta^b,
\eeq
where  $\eta_{ab}=\pm1\cdot\delta_{ab}$. 
The inverse metric tensor 
$g^{-1}=\sum_{a,b=1}^{N} \eta^{ab} X_{a}\otimes X_{b}$
(with  $\sum_{b} \eta^{ab}\eta_{bc}=\delta^{a}_{c}$ on the dual vector field basis such that $\theta^{b}(X_{a})=\delta^{b}_{a}$) allows to define an $\ca(G)$-bilinear contraction map $\tilde{g}:\Omega^{\otimes k}(G)\times\Omega^{\otimes s}(G)\to\Omega^{\otimes|k-s|}(G)$ given by 
\beq 
\label{cor}
\tilde{g}(\theta^{a_{1}}\otimes\ldots\otimes\theta^{a_{k}},\theta^{b_{1}}\otimes\ldots\otimes\theta^{ b_{s}})=
\left(\Pi_{i=1}^{k}\eta^{a_{i}b_{i}}\right)\theta^{b_{k+1}}\otimes\ldots\otimes\theta^{b_{s}}, 
\eeq
for $k<s$. In order to  set  such a map on the exterior algebra $\Omega(G)$,    we recall that  exterior forms are defined as the range of the classical antisymmetriser operators (those constructed from the  braiding given by the flip),   
\begin{align}
\theta^{a_{1}}\wedge\ldots\wedge\theta^{a_{s}}:=&A^{(s)}(\theta^{a_{1}}\otimes\ldots\otimes\theta^{a_{s}}) \nn \\
=&(A^{(k)}\otimes A^{(s-k)})\, B_{k,s-k}(\theta^{a_{1}}\otimes\ldots\otimes\theta^{a_{s}}),
\label{prAn}
\end{align} 
with $B_{k,s-k}:=\sum_{p_{j}\in\,S(k,s-k)}\mathrm{sign}(p_{j})\,p_{j}\,:\Omega^{\otimes s}(G)\to\Omega^{\otimes s}(G)$ the operator obtained in terms of the shuffles $S(k,s-k)\subset\,P(s)$, namely those   permutations $p_{j}$ of $s$ elements such that $p_{j}(1)<\ldots<p_{j}(k)$ and $p_{j}(k+1)<\ldots<p_{j}(s)$. The relation above proves that the expression
\beq
\label{Gexn}
\tilde{g}(\theta^{a_{1}}\wedge\ldots\wedge\theta^{a_{k}},\theta^{b_{1}}\wedge\ldots\wedge\theta^{ b_{s}})
:=\tilde{g}(A^{(k)}(\theta^{a_{1}}\otimes\ldots\otimes\theta^{a_{k}}),A^{(s)}(\theta^{b_{1}}\otimes\ldots\otimes\theta^{ b_{s}}))
\eeq
consistently defines a contraction $\tilde{g}:\Omega^{k}(G)\times\Omega^{s}(G)\to\Omega^{s-k}(G)$.
Fix an ordering $\vartheta=\theta^{1}\otimes\ldots\otimes\theta^{N}$ of the basis so to have the volume $N$-form $\mu=\mu^*:=A^{(N)}(\vartheta)$ (an orientation).  The corresponding Hodge duality $\star:\Omega^{k}(G)\to\Omega^{N-k}(G)$ is the operator whose action is
\begin{align}
\star(\xi):=&\,\frac{1}{k!}\,\tilde{g}(\xi,\mu) \nn \\
=&\,A^{(N-k)}(\tilde{g}(\xi,B_{k,N-k}\vartheta)).
\label{quid}
\end{align}
 It satisfies the identity
\beq
\label{quhs}
\star^{2}(\xi)=\mathrm{sgn}(g)(-1)^{k(N-k)}\xi
\eeq
on any $\xi\in\,\Omega^k(G)$. Here $\mathrm{sgn}(g)=\det(\eta_{ab})$ is the signature of the metric.

If the metric is in addition left-invariant, that is the 1-forms $\theta^{a}$ in \eqref{gmf} are left-invariant, the volume $\mu$ defines a left Haar measure on the group manifold, which allows to introduce a (left-invariant) sesquilinear inner product $\hs{~}{~}:\Omega^{k}(G)\times\Omega^{k}(G)\to\IC$ by
$\hs{\zeta}{\xi}:=\int_{\mu}\tilde{g}(\zeta^*,\xi);$ the Hodge duality associated to a left-invariant metric can  then be equivalently defined by the implicit equation
\beq
\hs{\zeta}{\xi}:=\int_{\mu}\zeta^*\wedge(\star\xi)
\label{Hsed}
\eeq
for any pair of $k$-forms $\zeta,\xi$.

\bigskip

The  problem of defining a Hodge duality operator on  exterior algebras associated to first order covariant differential calculi \`a la Woronowicz  over quantum groups has been studied following two different approaches, which parallel the two classical definition described above. 

A very general path is developed in \cite{hec99,hec03,hec00}. Given  a $*$-Hopf algebra $\ch$ equipped with a suitable bicovariant $N$-dimensional calculus, the theory of bicovariant bimodules allows, using the corresponding canonical braiding $\sigma$ (and its inverse $\sigma^{-}$), to introduce  a real $\sigma$-metrics (which generalises, as an homomorphism of $\ch$-bimodules, the  concept of real metric tensor on a classical group)  and a natural hermitian volume form. Using the (order $k$) antisymmetriser operators $\mathfrak{A}^{(k)}_{\pm}$ associated to $\sigma^{\pm}$, such metric structures extend to  contractions  between  higher order forms in the exterior algebras $\Omega_{\pm}(\ch)$;  Hodge duality operators $\star_{L,R}^{\pm}:\Omega^{k}_{\pm}(\ch)\to\Omega_{\pm}^{N-k}(\ch)$ are then given by (left, resp. right) $\ch$-linearly contracting
an exterior form with the volume. The quantum analogue of the classical \eqref{quhs} reads in this formulation
\beq 
\label{Idf}
\star_{L,R}^{\pm}\star_{L,R}^{\mp}(\xi)=\xi
\eeq 
on any  $\xi\in\,\Omega_{\pm}(\ch)$. 

A notion of left-invariant sesquilinear inner product on a large class of left-covariant exterior algebras $\Omega(\ch)$ over a $*$-Hopf algebra $\ch$ with an Haar state is introduced in \cite{kmt}, together with a natural definition of Hodge operator. This paper shows how specific constraints among inner products on $k$- and $(N-k)$-forms give an Hodge operator whose square has the classical (with $\mathrm{sgn}(g)=1$) spectrum \eqref{quhs}.   
This formulation has been also extended in \cite{ag10,ale09} to study Hodge operators on the left-covariant $\Omega(\sq)$ with $\sq$ the standard Podle\'s sphere.

\bigskip

The aim of the present  paper is to introduce  novel Hodge duality structures on the exterior algebra coming from the $4D_{+}$ first order bicovariant calculus \`a la Woronowicz on $\SU$, so to evolve the analysis started in \cite[\S3]{ag10}, from which  we get the geometrical setting and the notations.  
The quantum group $\SU$ and its exterior algebras $\Omega_{\pm}(\SU)$ are presented in the appendix \S\ref{eSU}: the  nontrivial spectral resolution of the braiding $\sigma^{\pm}$ reflects into a nontrivial spectral resolution of the corresponding antisymmetriser
operators $A^{(k)}_{\pm}$ (i.e. $A_{\pm}^{(k)}(\xi)=\lambda^{\pm}_{\xi}\,\xi$, with $\lambda_{\xi}^{\pm}\,\to\,k!$ in the formal classical limit $q\to1$),  whose actions do not moreover commute with the $*$-conjugation. 

This suggests to consider, in section \S\ref{s:Ho},  sesquilinear (instead of bilinear) contraction maps $\Gamma_{\pm}$ on the vector spaces $\Omega_{\pm {\rm inv}}(\SU)$ of left-invariant forms and to define, explicitly depending on  the spectra of the antisymmetrisers,  Hodge operators $\mT^{\pm},\,:\Omega^{k}_{\pm {\rm inv}}(\SU)\to\Omega^{4-k}_{\pm {\rm inv}}(\SU)$  which generalise  the second line out of \eqref{quid} to  the quantum setting. 

We restrict our analysis to contractions having a $\U(1)$-coinvariance, and we relate the   hermitianity and the reality of the contraction (i.e. a compatibility with the braiding and $*$-conjugation) to the behaviour of the corresponding Hodge operators.  This approach parallels \cite{ag10}, and somehow reverses the line described in \cite{hec99}, where the properties of a metric tensor give sufficient conditions for the Hodge duality to satisfy \eqref{Idf}. We shall here indeed study how natural requirements on the spectra  of the Hodge operators $\mT^{\pm}$ will constraint the family of contractions $\Gamma_{\pm}$ from which they are defined.

A contraction $\Gamma_{\pm}$ will  be called maximally hermitian provided the spectra of the square of the  corresponding Hodge operators have the same degeneracy of the antisymmetrisers; 
the main result of this analysis is that the Hodge operators $\mT^{\pm}$ corresponding to  a real  maximally hermitian contraction $\Gamma$ satisfy the identities, on the above diagonal basis for the antisymmetrisers, 
\begin{align}
&(\mT^{\pm})^{2}(\xi)=(-1)^{k(4-k)}(\mathrm{sgn}\,\Gamma)\,\frac{\lambda_{\xi}^{\pm}}{\lambda_{\xi^*}^{\pm}}\,\xi, \nn \\
&\mT^{\pm}\,\mT^{\mp}(\xi)=(-1)^{k(4-k)}(\mathrm{sgn}\,\Gamma)\,\xi
\label{idiT}
\end{align}
with a suitable definition of the signature of $\Gamma$: this induces to define Hodge dualities  $\star^{L,R}_{\pm}$  as the extension of the  $\mT^{\pm}$ (\emph{only} for real and  maximally hermitian contractions) to the whole exterior algebras $\Omega_{\pm}(\SU)$ via the requirement of left (or right) $\ASU$-linearity. 

From real maximally hermitian contractions we set (left-invariant) metric structures $\tilde{g}_{q}\in\,\mathfrak{G}$, whose signatures turn out to be constant:  
this is the same  happening to the class of  bicovariant real $\sigma$-metrics introduced on $\SU$ in \cite{hec99,hec03};  a comparison between the two notions of metrics  will reveal even more interesting similarities: this class of Hodge dualities allows furthermore to recover the quantum Casimir of $\su$ (the envelopping algebra dual to $\ASU$) as a Laplacian on $\SU$. 

We close this overview by saying that, to make the paper self-contained, appendix \S\ref{ass:a1} recalls basic notions on the theory of differential calculi and exterior algebras on  quantum groups.

\subsection*{Conventions and notations} 
We define a $q$-number as 
\beq
\label{qn}
[u]:=\frac{q^u-q^{-u}}{q-q^{-1}}
\eeq
 for any $q\neq 1$ and $u\in\,\IR$; for a coproduct
$\Delta$ we use the conventional Sweedler notation
$\Delta(x)=x_{(1)}\otimes x_{(2)}$ with implicit summation, iterated to  $(\id\otimes\Delta)\circ\Delta(x) = (\Delta\otimes\id)\circ\Delta(x) = x_{(1)}\otimes x_{(2)}\otimes x_{(3)}$, and so on.

\section{Hodge duality operators on $\SU$}
\label{s:Ho}

Motivated by the formulation we shall introduce in the quantum setting, we start by describing here  how Hodge dualities can be defined in the classical setting in terms of sesquilinear contraction maps.  On a $N$-dimensional real Lie group we call contraction an  $\ca(G)$-sesquilinear map $\Gamma:\Omega^1(G)\times\Omega^1(G)\to\ca(G)$ with
$\Gamma(f\,\phi,\phi^{\prime})=f^*\Gamma(\phi,\phi^{\prime})$ and $\Gamma(\phi,\phi^{\prime}\,f )=\Gamma(\phi,\phi^{\prime}) f$ for  
$f\in\,\ca(G)$, which we extend  (see \eqref{cor},\eqref{prAn}) to  
$\Gamma:\Omega^{k}(G)\times\Omega^{s}(G)\to\Omega^{s-k}(G)$ via
\beq
\label{Gex}
\Gamma(\phi^{a_{1}}\wedge\ldots\wedge\phi^{a_{k}},\phi^{b_{1}}\wedge\ldots\phi^{ b_{s}})
:=\Gamma(A^{(k)}(\phi^{a_{1}}\otimes\ldots\otimes\phi^{a_{k}}),A^{(s)}(\phi^{b_{1}}\otimes\ldots\phi^{ b_{s}})).
\eeq
Consider an ordered left-invariant basis  $\{\phi^1,\ldots,\phi^N\}\subset\Omega^1(G)$, and define $\vartheta=m\,\phi^1\otimes\ldots\otimes\phi^{N}$ with $m\,\in\IC$ so that $\mu=\mu^*=\mathfrak{A}^{(N)}(\vartheta)$ is an hermitian volume form for the calculus. Given $\xi\in\,\Omega^{k}(G)$, the expression 
\beq
\label{deft}
T(\xi):=(A^{(N-k)}(\Gamma(\xi,B_{k,N-k}\theta)))^*
\eeq
defines an $\ca(G)$-linear operator $\Omega^{k}(G)\to\Omega^{N-k}(G)$.
In order to recover the operator $T$ as an Hodge operator corresponding to a metric structure, define the tensor $\tilde{g}:\Omega^{\otimes 2}(G)\to\ca(G)$ on the above basis as
\beq
\label{degti}
\tilde{g}(\phi^{a},\phi^{b}):=\Gamma(\phi^{a*},\phi^b);
\eeq
it is then immediate to see that 
 $T(1)=\mu^*$
and $T(\mu^*)=(m^{2}\det\,\tilde{g})^*$, where $\det\,\tilde{g}$ is the determinant of the matrix $\tilde{g}(\phi^{a},\phi^b)$. 

Imposing on the sesquilinear $\Gamma$ both an hermitianity and a reality condition  will impose, as in \cite{ag10}, both a symmetry and a reality condition on the tensor $\tilde{g}$:
\begin{itemize}
\item 
The sesquilinear contraction  map $\Gamma$ is said hermitian provided the operator \eqref{deft}
satisfies
\beq
\label{dHe}
T^2(\phi)=(-1)^{N-1}\,T^2(1)\,\phi
\eeq
for any $\phi\in\,\Omega^1(G)$. Easy calculations using the properties of the antisymmetriser (and shuffle) operators prove that 
\beq
\label{eqhs}
T^2(\phi)=(-1)^{N-1}\,\{T^2(1)\}\,\phi\qquad\Leftrightarrow\qquad\tilde{g}(\phi^a,\phi^b)=\tilde{g}(\phi^{b},\phi^a):
\eeq
hermitian sesquilinear contractions $\Gamma$ corresponds to symmetric bilinear tensors $\tilde{g}$. 
\item  
The sesquilinear contraction map $\Gamma$ is said real provided 
\beq
\label{crG}
T(\phi^*)=(T(\phi))^*
\eeq
for  any $\phi\in\,\Omega^1(G)$. Such a reality conditions on the sesquilinear map $\Gamma$ is equivalent to a reality condition for the bilinear tensor $\tilde{g}$, namely 
\beq
\label{eGg}
\tilde{g}^*(\phi^a,\phi^b)=\tilde{g}(\phi^{b*},\phi^{a*}),
\eeq
and moreover implies  that the action of the operator $T$ will commute with the $*$-conjugation, i.e. $T(\xi^*)=(T(\xi))^*$,  on any  $\xi\in\Omega(G)$.  
\end{itemize}
Upon fixing the scale parameter by  $m^2\det\,\tilde{g}=\mathrm{sgn}(\tilde{g})$, it is  straightforward to see that the operator $T$ associated to hermitian and real contraction maps $\Gamma$ is the Hodge operator  on $\Omega(G)$ corresponding to   the inverse of the real metric tensor $\tilde{g}$ given by $\Gamma$ via \eqref{degti}.

\begin{rema}
We recall, since we shall use it in the following, that the action of the operator $T$ coincides with the action of the operator $L:\Omega^{k}(G)\to\Omega^{N-k}(G)$ defined in the same setting  by 
\beq
\label{Ldf}
L(\xi):=\frac{1}{k!}\,\Gamma^*(\xi,\mu)
\eeq
on $\xi\in\,\Omega^{k}(G)$; hermitianity and reality  condition are given by:
\begin{align}
&\phi^*\wedge L(\phi^{\prime})=\Gamma(\phi,\phi^{\prime})\mu,  \nn \\
&L(\phi^*)=(L(\phi))^*
\label{hrL}
\end{align}
on any left-invariant 1-form $\phi$. 
\end{rema}

\subsection{(Sesquilinear) Contraction maps in the quantum setting}
\label{ss:cm}
In order to introduce Hodge operators on the exterior algebras $\Omega_{\pm}(\SU)$ generated by the antisymmetriser operators $A^{(k)}_{\pm}$ (as described in appendix \S\ref{se:4dc})    we shall start by considering left $\ASU$-invariant contractions as   sesquilinear maps 
$$
\Gamma:\Omega^1_{{\rm inv}}(\SU)\times\Omega^1_{{\rm inv}}(\SU)\to \IC,
$$
 with $\Gamma(\lambda^*\,\omega,\omega^{\prime})=\Gamma(\omega,\lambda\,\omega^{\prime})=\lambda\,\Gamma(\omega,\omega^{\prime})$ for any $\lambda\in\,\IC$ and $\omega,\omega^{\prime}\in\,\Omega^{1}_{{\rm inv}}(\SU)$, which we  consistently extend to contractions   
 $$
 \Gamma_{\pm}:\Omega^{ k}_{\pm {\rm inv}}(\SU)\times\Omega^{N}_{\pm {\rm inv}}(\SU)\to\Omega^{ N-k}_{\pm {\rm inv}}(\SU),
 $$
  where the  shuffle operators are given by $B^{\pm}_{k,N-k}:=\sum_{p_{j}\in\,S(k,N-k)}\mathrm{sign}(p_{j})\mathfrak{P^{\pm}}(p_{j})$ in terms of permutation operators $\mathfrak{P}^{\pm}(p)$ defined by the braiding $\sigma$ of the bicovariant calculus on $\SU$ or of its inverse $\sigma^{-1}$ (see \S\ref{ass:a1}). 

From \eqref{om4} we set 
$\Omega^{\otimes 4}_{{\rm inv}}(\SU)\,\ni\,\vartheta:=im\,\omega_{-}\otimes\omega_{+}\otimes\omega_{0}\otimes\omega_{z}$ with $m\neq0\in\,\IR$, so to have, from \eqref{sp3-},
$\mu\,=\,\mu^*\,=\,
\check{\mu}\,=\,\check{\mu}^*$. 
On a left-invariant basis in $\Omega_{\pm}(\SU)$ which diagonalises (appendix \S\ref{se:4dc}) the action of the antisymmetrisers,
\beq
\label{eiA}
A_{\pm}^{(k)}(\xi)=\lambda^{\pm}_{\xi}\,\xi,
\eeq 
we define the operator 
$\mT^{\pm}:\Omega_{\pm {\rm inv}}^{k}(\SU)\to\Omega_{\pm {\rm inv}}^{N-k}(\SU)$ by 
\begin{align}
& \mT^{\pm}(\xi):=\frac{\lambda^{\pm}_{\xi^*}}{\lambda^{\pm}_{\xi}}\,(A^{(N-k)}_{\pm}(\Gamma_{\pm}(\xi,B^{\pm}_{k,N-k}\vartheta)))^*.  
\label{dTL}
\end{align}
These expression generalises to the quantum setting  the classical \eqref{deft} taking into account the specific spectrum of the quantum antisymmetrisers. 

Since we are interested in Hodge dualities which will give Laplacians  whose action restricts \eqref{clnm} to any $\cl_{n}\subset\ASU$, throughout the paper we shall analyse the properties of the operators defined in \eqref{dTL}, corresponding to the  class of left-invariant contractions which are  (right) $\U(1)$-coinvariant. 
Given $\delta_{R}^{(1)}\omega_{a}=\omega_{a}\otimes z^{n_{a}}$ with $\delta_{R}^{(1)}:\Omega^1(\SU)\to\Omega^1(\SU)\otimes\ca(\U(1))$ the extension to 1-forms of the $\U(1)$-coaction \eqref{cancoa}, a contraction map is said $\U(1)$-coinvariant provided
\beq
\label{dU}
n_{a}\neq n_{b}\qquad
\qquad\Rightarrow\qquad\qquad\Gamma_{ab}=\Gamma(\omega_{a},\omega_{b})=0. 
\eeq
The coefficient of a $\U(1)$-coinvariant contraction map can be written as (in the ordering defined above):
\beq
\label{Ggc}
\Gamma_{ab}=\left(\begin{array}{cccc}
\alpha & 0 & 0 & 0 \\ 0 & \beta & 0 & 0 \\ 0 & 0 & \nu & \epsilon \\ 0 & 0 & \xi & \gamma
\end{array}\right)
\eeq
where all the parameters are complex numbers.

\subsection{The Hodge operator $\mT^{+}$}
\label{ss:T}
We start a more detailed  analysis of the Hodge operator $\mT^+$.  From the definition \eqref{dTL} one has
\begin{align}
&\mT^+(1)=\mu, \nn \\
&\mT^+(\mu)=\Gamma_{+}^*(\mu,\theta)=m^2(\alpha\,\beta)^*(\nu\,\gamma-\epsilon\,\xi)^*
\label{Tfa}
\end{align}
while, on the basis of left-invariant 1-forms,
\beq
\label{T1}
\mT^+\left(\begin{array}{c} \omega_{-} \\ \omega_{+} \\ \omega_{0} \\ \omega_{z} \end{array}\right)\,=\,i\,m
\left(\begin{array}{cccc}
0 & \alpha^*  & 0 & 0 \\ -\beta^* & 0 & 0 & 0  \\ 0 & 0 & -\nu^* & \epsilon^* \\ 0 & 0 & -\xi^* & \gamma^*
\end{array}\right)
\left(\begin{array}{c} \chi_{-} \\ \chi_{+} \\ \chi_{0} \\ \chi_{z} \end{array} \right).
\eeq
As a reality condition, we set the classical relation \eqref{crG}, which  imposes specific constraints on the contraction coefficient  $\Gamma_{+}$ \eqref{Ggc}:
\beq
\begin{array}{c} \mT^+(\omega_{a}^*)=(\mT^+(\omega_{a}))^* \\ \end{array} \qquad\Leftrightarrow\qquad \begin{array}{l} \beta^*=q^2\alpha, \\
(\nu,\epsilon,\xi,\gamma)\,\in\,\IR. \end{array}
\label{reT}
\eeq
For a real  contraction map the action of the Hodge operator $\mT^+$ on 3-forms  is
\begin{align}
&\mT^+(\chi_{-})=-im\,q^{-2}\beta^*(\nu\,\gamma-\epsilon\,\xi)\,\omega_{+}, \nn \\
&\mT^+(\chi_{+})=im\,q^2\alpha^*(\nu\,\gamma-\epsilon\,\xi)\,\omega_{-}, \nn \\
&\mT^+\left(\begin{array}{c}\chi_{0} \\ \chi_{z} \end{array} \right)=im\left(\begin{array}{cc}
\gamma^2\nu+(\alpha\,\beta)^*((1-q^2-q^{-2})\gamma-2(q+q^{-1})^2\nu-(q^2-q^{-2})\xi) & (\alpha\,\beta)^*\xi \\  \epsilon^2\xi+(\alpha\,\beta)^*((q^2-q^{-2})\nu+(1-q^2-q^{-2})\epsilon) & (\alpha\,\beta)^*\nu 
\end{array}
\right) \left(\begin{array}{c} \omega_{0} \\ \omega_{z} \end{array}\right).
\label{T3}
\end{align} 
This allows to set the hermitianity conditions,
\beq
\label{hT+}
\begin{array}{c} (\mT^{+})^{2}(\omega_{a})=-(\mT^+)^2(1)\,\omega_{a}\end{array} \quad\Leftrightarrow\quad\begin{array}{l} \beta=q^2\alpha, \\ \xi=\epsilon, \\ \epsilon^3+\alpha\,\beta((q^2-q^{-2})\nu-(q-q^{-1})^2\epsilon)=0, \\ \gamma^2\nu=\alpha\,\beta((q^2-q^{-2})\epsilon-2(q+q^{-1})^2\nu-(q-q^{-1})^2\gamma), \end{array} 
\eeq
which, together with the reality conditions above, give  three families of real and hermitian $\U(1)$-covariant contractions which are invertible (that is, $\mT^+(\mu)\neq0$):
\begin{enumerate}[a)]
\item   
the first one is 
\beq
\label{as}
\begin{array}{l} (\alpha\neq0,\; \beta=q^2\alpha,\; \nu=0, \; \gamma\neq0,\; \epsilon\neq0, \;\xi=\epsilon )\,\in\,\IR, \\ \epsilon^2=(q-q^{-1})^2\alpha\,\beta, \\ (q+q^{-1})\epsilon + (q-q^{-1})\gamma=0;
\end{array}
\eeq
\vskip3pt
\item the second one is:
\beq
\label{bs}
(\alpha\neq0,\; \beta=q^2\alpha,\; \nu\neq0, \; \gamma\neq0,\; \epsilon\neq0, \;\xi=\epsilon )\,\in\,\IR
\eeq
where the third and fourth out of \eqref{hT+} are left implicit. These can be solved if $(q-q^{-1})^2\alpha\,\beta>\epsilon^2$;
\vskip 3pt
\item the third one is: 
\beq
\label{cs}
\begin{array}{l} (\alpha\neq0,\; \beta=q^2\alpha,\; \nu\neq0, \; \gamma=0,\; \epsilon\neq0, \;\xi=\epsilon )\,\in\,\IR, \\ 2\,\epsilon^2=3(q-q^{-1})^2\alpha\,\beta, \\ 2(q+q^{-1})\nu + (q-q^{-1})\epsilon=0.
\end{array}
\eeq
\end{enumerate}
We notice that, since the braiding of the calculus is not the classical flip, hermitianity and reality constraints are not  only expressed  by linear relations among the coefficients of the contractions. 

The action of the Hodge operator $\mT^+$ on left-invariant 2-forms is:
\beq
\label{T2f}
\mT^+\left(\begin{array}{c} \varphi_{+} \\ \kappa_{-} \end{array} \right)= -im\,\alpha
\left(\begin{array}{cc}  q^2\epsilon+(\frac{q^2+1-q^{-2}}{q^2-1})\nu & (\frac{1}{1-q^2})\nu \\ 
(q^4+q^2)\epsilon+(q^4-q^2)\gamma+(\frac{q^6-q^2+1}{q^2-1})\nu & -\epsilon +(\frac{1}{1-q^2})\nu
 \end{array}\right)\left(\begin{array}{c} \varphi_{+} \\ \kappa_{-} \end{array}\right),
\eeq
\beq
\label{T2p}
\mT^+\left(\begin{array}{c} \varphi_{-} \\ \kappa_{+} \end{array} \right)= -im\,\beta
\left(\begin{array}{cc}  -q^{-2}\epsilon+(\frac{q^{-2}+1-q^{2}}{q^{-2}-1})\nu & (\frac{1}{1-q^{-2}})\nu \\ 
-(q^{-4}+q^{-2})\epsilon+(q^{-4}-q^{-2})\gamma+(\frac{q^{-6}-q^{-2}+1}{q^{-2}-1})\nu & \epsilon +(\frac{1}{1-q^{-2}})\nu
 \end{array}\right)\left(\begin{array}{c} \varphi_{-} \\ \kappa_{+} \end{array}\right), 
\eeq
\begin{align}
&\mT^+(\psi_{-})=-\frac{im}{1+q^2}\left(\{(\frac{1-q^2-q^{-2}}{1-q^{-2}})\nu\,\gamma+(\frac{2-q^2}{1-q^{-2}})\epsilon^2+(q^2-1)(q^2+q^{-2}-1)\alpha\,\beta\}\psi_{-}\right. \nn
\\ &\qquad\qquad\qquad\qquad\qquad+\left.\{(\frac{2q^2-1}{1-q^{-2}})\nu\,\gamma+(\frac{-q^4+q^2-1}{1-q^{-2}})\epsilon^2+(q^2-1)(q^4+1-q^2)\alpha\,\beta\}\psi_{+}
\right) \nn \\
&\mT^+(\psi_{+})=-\frac{im}{1+q^2}\left(\{(\frac{1-2q^{-2}}{1-q^{-2}})\nu\,\gamma+(\frac{1-q^{-2}-q^{-4}}{1-q^{-2}})\epsilon^2+(q^{-2}-1)(q^2+q^{-2}-1)\alpha\,\beta\}\psi_{-} \right.\nn \\ 
&\qquad\qquad\qquad\qquad\qquad+\left.\{(\frac{q^2-1+q^{-2}}{1-q^{-2}})\nu\,\gamma+(\frac{q^{-2}-2}{1-q^{-2}})\epsilon^2+(1-q^2)(q^2+q^{-2}-1)\alpha\beta\}\psi_{+}\right).
\label{T2s}
\end{align}
Every family  of real and hermitian contractions \eqref{as}-\eqref{cs} above will clearly give a specific spectrum to the restriction of Hodge operator $(\mT^+)^2$ on $\Omega^2_{+inv}$. Among them, we look for those conditions giving such spectra a specific degeneracy, thus recalling the classical case \eqref{quhs}, when the spectrum of square of a Hodge operator is totally degenerate. 

We define then a contraction $\Gamma_{+}$ \emph{real} and \emph{maximally hermitian} (with respect to the Hodge operators $\mT^+$) provided:
\begin{itemize}
\item it is real and hermitian;
\item the operator $(\mT^+)^2$ has the same degeneracy of the antisymmetrisers  $A^{(k)}_{+}$, namely the square of the Hodge operator $\mT^+$ is constant on each eigenspace of the antisymmetrisers $A^{(k)}_{+}$. 
\end{itemize}
Direct calculations show that the $\U(1)$-coinvariant contraction  in \eqref{Ggc}  is real and maximally hermitian provided the constraints \eqref{as} -- given by the set of  conditions a) above -- are fullfilled. The action of the corresponding Hodge operator $\mT^+$ turns out to be,  from the expressions \eqref{T2f}-\eqref{T2s}, as:
\begin{align}
&\mT^+(\varphi_{+})=-im\,q^2\alpha\,\epsilon\,\varphi_{+}, \nn \\
&\mT^+(\kappa_{+})=-im\,q^2\alpha\,\epsilon\,\kappa_{+}, \nn \\
&\mT^+(\psi_{+})=im(q^2-1)\alpha\,\beta\,\psi_{+}, \nn \\
&\mT^+(\varphi_{-})=im\,\alpha\,\epsilon\,\varphi_{-}, \nn \\
&\mT^+(\kappa_{-})=im\,\alpha\,\epsilon\,\kappa_{-}, \nn \\
&\mT^+(\psi_{-})=-im(1-q^{-2})\alpha\,\beta\,\psi_{-}.
\label{T2ms}
\end{align}
Given this class of Hodge operators, the reality condition \eqref{reT} can be extended to a more general compatibility 
of the action of  $\mT^+$  with the $*$-conjugation in  $\Omega_{+{\rm inv}}(\SU)$; for any left-invariant "eigenform" \eqref{eiA} of the antisymmetriser $A_{+}^{(k)}$ one has 
\beq
\lambda^{+}_{\xi^*}\,(\mT^+(\xi))^*=\lambda_{\xi}^{+}\mT^+(\xi^*).
\label{rTt}
\eeq 
We close the analysis of the Hodge operator $\mT^+$ by noticing that, on the same  basis, for a real and maximally hermitian contraction one has
\beq
(\mT^+)^2(\xi)=(-1)^{k(4-k)}\frac{\lambda^{+}_{\xi}}{\lambda^{+}_{\xi^*}}\,(\mT^+)^2(1)\,\xi.
\label{pq}
\eeq
In order to clarify this relation we can define  the quantum determinant of the   contraction map $\Gamma_{+}$ as
\beq
\mathrm{det}_{q}\,\Gamma_{+}:=\Gamma_{+}(i\,\omega_{-}\otimes\omega_{+}\otimes\omega_{0}\otimes\omega_{z},
i\,\omega_{-}\wedge\omega_{+}\wedge\omega_{0}\wedge\omega_{z})
\label{qde}
\eeq
and fix (up to a sign) the scale parameter $m\in\,\IR$ imposing 
\beq
(\mT^+)^2(1)=\mathrm{sgn}\,\Gamma_{+}\qquad\Leftrightarrow\qquad
m^2(\alpha\,\beta\,\epsilon^2)=1,
\label{boh}
\eeq
with $\mathrm{sgn}\,(\Gamma_{+})=(\mathrm{det}_{q}\Gamma_{+})\backslash\left|\mathrm{det}_{q}\Gamma_{+}\right|$. This choice allows to write the relation \eqref{pq} as
\beq
\label{sTg}
(\mT^+)^2(\xi)=(-1)^{k(4-k)}(\mathrm{sgn}\,\Gamma_{+})\frac{\lambda^{+}_{\xi}}{\lambda^{+}_{\xi^*}}\,\xi,
\eeq
which generalises the classical \eqref{quhs} to the quantum setting, where the braiding $\sigma$ associated to the calculus on $\SU$ has the non trivial spectral decomposition given in \eqref{spb}.

\subsection{The Hodge operator $\mT^{-}$}
It is now easy to study the Hodge operator $\mT^-$ on $\Omega_{-{\rm inv}}(\SU)$,  following  the path outlined for $\mT^+$. Given the $\U(1)$-coinvariant contraction map $\Gamma$ in \eqref{Ggc}, from  the second relation out of \eqref{sp3-} we have
\begin{align}
&\mT^-(1)=\check{\mu}=\mu=\mT^+(1), \nn  \\ 
&\mT^-(\check{\mu})=\mT^+(\mu),
\label{mT04}
\end{align}
while, from the first relation out of \eqref{sp3-},
\beq
\mT^-(\omega_{a})=\mT^+(\omega_{a}).
\label{mT1}
\eeq
Some algebra allows to work out the action of $\mT^-$ on $\Omega^3_{-{\rm inv}}(\SU)$, so to obtain the conditions of reality and hermitianity of the contraction $\Gamma_{-}$; for real and hermitian contraction maps one has (with $\chi_{a}=\check{\chi}_{a}$, recall \eqref{sp3-})
\beq
\label{e3mp}
\mT^-(\chi_{a})=\mT^+(\chi_{a}).
\eeq
Some more algebra gives then the action of $\mT^-$ on  $\Omega_{{\rm -inv}}^2$ and proves  indeed  that  the family of real and maximally hermitian contractions $\Gamma_{-}$ for $\mT^-$ coincides with that obtained above by the Hodge operator $\mT^+$. For those contractions:
\begin{align}
&\mT^-(\check\varphi_{+})=-im\,\alpha\,\epsilon\,\check\varphi_{+}, \nn \\
&\mT^-(\check\kappa_{+})=-im\,\alpha\,\epsilon\,\check\kappa_{+}, \nn \\
&\mT^-(\check\psi_{+})=-im(q^{-2}-1)\alpha\,\beta\,\check\psi_{+}, \nn \\
&\mT^-(\check\varphi_{-})=im\,q^2\alpha\,\epsilon\,\check\varphi_{-}, \nn \\
&\mT^-(\check\kappa_{-})=im\,q^2\alpha\,\epsilon\,\check\kappa_{-}, \nn \\
&\mT^-(\check\psi_{-})=im(1-q^{2})\alpha\,\beta\,\check\psi_{-}.
\label{T2mms}
\end{align}
The Hodge operator $\mT^-$ corresponding to real and maximally hermitian contractions  satisfies, on 
the left-invariant basis of $\Omega_{-}(\SU)$ given by \eqref{eiA}, the same reality condition \eqref{rTt} that $\mT^+$ satisfies, namely 
\beq
\lambda^{-}_{\xi^*}\,(\mT^-(\xi))^*=\lambda_{\xi}^{-}\mT^-(\xi^*),
\label{rTm}
\eeq 
and its square gives, in complete analogy  to \eqref{pq},
\beq
(\mT^-)^2(\xi)=(-1)^{k(4-k)}\frac{\lambda^{-}_{\xi}}{\lambda^{-}_{\xi^*}}\,(\mT^-)^2(1)\,\xi.
\label{pq-}
\eeq
It is now immediate to check that, with 
$$
\mathrm{det}_{q}\,\Gamma_{-}:=\Gamma_{-}(i\,\omega_{-}\otimes\omega_{+}\otimes\omega_{0}\otimes\omega_{z},
i\,\omega_{-}\wedge\omega_{+}\wedge\omega_{0}\wedge\omega_{z})=\mathrm{det}_{q}\,\Gamma_{+}
$$ 
(where the last equality comes from \eqref{mT04}), from  the condition 
$(\mT^-)^2(1)=\mathrm{sgn}\,\Gamma_{-}\;\Leftrightarrow\;
m^2(\alpha\,\beta\,\epsilon^2)=1$
one can write
\beq
\label{sTgm}
(\mT^-)^2(\xi)=(-1)^{k(4-k)}(\mathrm{sgn}\,\Gamma_{-})\frac{\lambda^{-}_{\xi}}{\lambda^{-}_{\xi^*}}\,\xi.
\eeq
Although $\U(1)$-coinvariant real and maximally hermitian contractions $\Gamma_{\pm}$ \emph{coincide},  the actions (as relations \eqref{T2ms} and \eqref{T2mms} show) on $\Omega_{\pm{\rm  inv}}(\SU)$ of the corresponding Hodge operators $\mT^{\pm}$ \emph{do not}.  The $\Omega_{- {\rm inv}}(\SU)\simeq\Omega_{+{\rm inv}}(\SU)$ isomorphism described in \eqref{spe2-} and \eqref{sp3-} indeed allows to recover  that, with the obvious position $\mathrm{det}_{q}\,\Gamma_{+}=\mathrm{det}_{q}\,\Gamma_{-}$ for this class of contractions,
\beq
\mT^+\mT^-(\xi)=\mT^-\mT^+(\xi)=(-1)^{k(4-k)}(\mathrm{sgn}\,\Gamma)\,\xi;
\label{fq}
\eeq
this relation shows, in a quantum setting with $\sigma^2\neq1$,  the closest similarity to the classical \eqref{quhs}.

The last step is to extend the Hodge operators to the whole exterior algebras $\Omega_{\pm}(\SU)$ by  requiring them to be left (resp. right) $\ca(\SU)$-linear. We  define then \emph{Hodge duality operators} 
$\star_{\pm}^{L,R}:\Omega^{k}(\SU)\to\Omega^{4-k}(\SU)$ by
\begin{align}
&\star^{L}_{\pm}(x\,\omega):=x\,\mT^{\pm}(\omega), \nn \\
&\star^{R}_{\pm}(\omega\,x):=\mT^{\pm}(\omega)\,x ,
\label{shc}
\end{align}
(with $x\in\,\ASU$ and $\omega\in\,\Omega_{{\rm inv}}(\SU)$) 
for the class of $\U(1)$-coinvariant  real and maximally hermitian contraction map $\Gamma$ 
considered above. We summarise the results of this section as a proposition.
\begin{prop}
Given the Hodge operators $\mT^{\pm}$ defined in \eqref{dTL}, their corresponding sets of real and maximally hermitian contractions $\Gamma_{\pm}$ do coincide. It is always possible to rescale (with $m\in\,\IR$) the volume form $\mu=\check{\mu}=im\,\omega_{-}\wedge\omega_{+}\wedge\omega_{0}\wedge\omega_{z}$ such that relations \eqref{rTt}, \eqref{sTg} as well as \eqref{rTm}, \eqref{sTgm} hold. Their extension to the Hodge duality operators moreover satisfy
$$
\star^{L,R}_{\pm}\star^{L,R}_{\mp}(\phi)=(-1)^{k(4-k)}(\mathrm{sgn}\,\Gamma)\phi
$$
for any $\phi\in\,\Omega^{k}(\SU)$.
\end{prop}

\section{Symmetric contractions and Laplacians on $\SU$}
\label{senu}

For the class of real and maximally hermitian contractions  one has $\mathrm{det}_{q}\Gamma=-q^2(1-q^2)^2\alpha^4$, so that $\mathrm{sgn}\,\Gamma=-1$: we recover here a phenomenon which also the formulation in \cite{hec99,hec03} presents. 
For a deeper analysis of the relations between the two formulations, following  \eqref{degti} we define a  metric structure on the left-invariant part of the  first order differential  calculus over $\SU$ as a bilinear  map $\tilde{g}_{q}:\Omega^{\otimes 2}_{{\rm inv}}(\SU)\to\IC$ set by
\beq
\tilde{g}_{q}(\omega_{a},\omega_{b}):=\Gamma(\omega^*_{a},\omega_{b})
\label{gtq}
\eeq
 associated to  contraction maps $\Gamma$ above. 
Relations \eqref{as} show that this class of metrics, that we denote by $\mathfrak{G}$, can be parametrised by $(\IR\backslash {0})\times \IZ_{2}$, that is a nonzero real parameter together with the choice of a sign. In the matrix notation \eqref{Ggc},
\beq
\tilde{g}_{q\,ab}=\left(\begin{array}{cccc}  0 & q^2\ma & 0 & 0 \\ \ma & 0 & 0 & 0 \\ 0 & 0 & 0 &\pm(1-q^2)\ma \\ 0 & 0 & \pm(1-q^2)\ma & \pm(q^2+1)\ma \end{array}\right)
\label{cop}
\eeq
with $\ma=-\alpha\neq0$. 
In the formulation introduced and developed in \cite{hec99,hec03},  Hodge duality operators acting on     $\Omega(\SU)$ are defined from non-degenerate $\sigma$-metrics, which are maps $\mg:\Omega^{\otimes2}(\SU)\to\ASU$ satisfying the conditions
\begin{itemize}
\item $\mg$ is an homomorphism of the $\ASU$-bimodule $\Omega^{\otimes2}(\ASU)$;
\item $\mg$ is non-degenerate;
\item $\mg\circ\sigma^{\pm}=\mg$ (symmetry condition);
\item the equality $(\mg\otimes 1)\circ\sigma^{\pm}_{2}=(1\otimes\mg)\circ\sigma^{\mp}_{1}$ on $\Omega^{\otimes3}(\SU)$  holds;
\end{itemize}
$\sigma$-metrics are called real provided $(\mg(\phi,\phi^{\prime}))^*=\mg(\phi^{\prime*},\phi^*)$, and left-covariant provided $\Delta\circ\mg=(1\otimes\mg)\Delta_{L}^{(2)}.$

In order to compare these two definitions, we restrict our analysis to the set -- that we denote $\mathfrak{G}_{\sigma}$ -- of  real $\sigma$-metrics  given by $\mg(\omega_{a},\omega_{b})\in\,\IC$ with $\omega_{j}$ left-invariant 1-forms (such are left-covariant), and which are in addition right $\U(1)$-coinvariant, namely  (using the notation of \eqref{dU}) those satisfying $n_{a}+n_{b}\neq0\;\Rightarrow\mg(\omega_{a},\omega_{b})=0$ with $\delta_{R}^{(1)}(\omega_{a})=z^{n_{a}}\otimes\omega_{a}$.

An  explicit calculation proves that $\mathfrak{G}_{\sigma}\subset\mathfrak{G}$: any   $\mg\in\,\mathfrak{G}_{\sigma}$ coincides with  a  metric structure $\tilde{g}_{q}$ in \eqref{cop} corresponding to the contraction with $\epsilon=(1-q^2)\alpha,\;\gamma=(1+q^2)\alpha.$

\bigskip

We gain a further perspective on metrics in $\mathfrak{G}_{\sigma}$ and $\mathfrak{G}\backslash\mathfrak{G}_{\sigma}$ by focussing on the spectra of the Laplacians associated  to the  Hodge dualities
\eqref{shc}.
The equalities \eqref{mT1}, \eqref{e3mp}   allow to   define  operators $\Box^{L,R}:\ASU\to\ASU$ whose action can be written in terms of the basic derivations \eqref{Lq}, \eqref{d4},
\begin{align}
&\Box^{L}x:=-\star^{L}\,\dd\,\star^L\,\dd \,x\,=\alpha\{L_{+}L_{-}+q^2L_{-}L_{+}\mp(1+q^2)L_{z}L_{z}\pm2(q^2-1)L_{0}L_{z}\}\lt x, \nn \\
&\Box^{R}x:=-\star^{R}\,\dd\,\star^R\,\dd \,x\,=\alpha\{q^2R_{+}R_{-}+R_{-}R_{+}\mp(1+q^2)R_{z}R_{z}\pm2(q^2-1)R_{0}R_{z}\}\lt x.
\label{Ldc}
\end{align}
These operators are  diagonal on the vector space basis $\ASU\supset\cl_{n}\,\ni\,\phi_{n,J,l}=(c^{J-n/2}a^{*J+n/2})\rt E^l$ with $n\in\,\IZ,\,J=|n|/2+\IN,\,l=0,\ldots,2J$, and their spectrum is given in \cite[\S5]{ag10}. 
We report the analysis on  Laplacians $\Box^L$ in more details, being the pattern for $\Box^R$ analogue.

 For a metric in $\mathfrak{G}_{\sigma}$ we have, from \eqref{Lq} and the \eqref{casbis} which gives the quantum Casimir of $\su$,
\beq
\label{sur}
\Box^{L}x=\alpha\{L_{+}L_{-}+q^2L_{-}L_{+}+(1+q^2)L_{z}L_{z}-2(q^2-1)L_{0}L_{z}\}\lt x=2q\,\alpha\,L_{0}\lt x,
\eeq
with 
$$
L_{0}\lt \phi_{n,J,l}=[J][J+1]\phi_{n,J,l}.
$$
 The action of this Laplacian\footnote{The equality $\Box^{R}\,x=\Box^{L}\,x$ if $\tilde{g}_{q}\in\,\mathfrak{G}_{\sigma}$ is immediate to check.} reduces in the  classical limit  to the action of the quadratic Casimir $C$ of the Lie algebra $\mathfrak{g}=\mathfrak{su(2)}$. It is indeed well known that $C$  coincides with the Laplacian corresponding to the Cartan-Killing metric on $\mathrm{SU}(2)$: it seems then  meaningful to assume that a metric   $\tilde{g}_{q}\in\,\mathfrak{G}_{\sigma}$ reduces in the classical limit to a  Cartan-Killing metric on $\mathrm{SU}(2)$.   
Along the same path it is also immediate to recover that, although real and maximally hermitian contractions have the same signature in the quantum setting, corresponding metrics in $\mathfrak{G}_{\sigma}$ and $\mathfrak{G}\backslash\mathfrak{G}_{\sigma}$ have different signatures in the classical limit. The discrete spectrum of Laplacians corresponding to metrics in 
$\mathfrak{G}\backslash\mathfrak{G}_{\sigma}$ moreover ranges on the whole real line: although  seems it  meaningless to use them to model energy levels of a stable quantum system, this class of Laplacians show that the formalism developed here consistently introduce some features of a non Riemannian geometry on $\SU$.

\begin{rema}
The Laplacian \eqref{sur} shows an interesting relation with the Laplacian $\Delta_{q}$ introduced in \cite{bk96} on $\SU$ in terms of the $R$-matrix structure of $\ASU$ and without any reference to a differential calculus. If we fix $\alpha=(q^2+1)^{-1}$ in \eqref{sur}, then 
\beq
\label{coLa}
\Delta_{q}=\Box^{L}+\left(\frac{q-q^{-1}}{q+q^{-1}}\right)^2L_{0}^2 
\eeq
as elements in $\su$. To cast the analysis in \cite{bk96} in the perspective of a general study of Dirac operators on $\SU$ we refer to \cite{grev} and reference therein. 
\end{rema}




\subsection{A comparison}
\label{ss:Lh}
We describe now how the set of metrics $\mathfrak{G}$ can be characterised by an analysis of the Hodge operators  $\mL^{\pm}$ on $\Omega_{\pm inv}(\SU)$  defined in \cite{ag10} as 
\beq
\label{qLope}
\mL^{\pm}(\xi):=\frac{1}{\lambda^{\pm}_{\xi}}\,\Gamma_{\pm}^*(\xi,\mu).
\eeq
This epression clearly generalises the classical \eqref{Ldf}.
From the definition of shuffle operators one has $B^{\pm}_{1,N-1}=1$, so that (with $\mu=\check\mu$)
\begin{align}
&\mL^{\pm}(1)=\mu, \nn \\
&\mL^{\pm}(\mu)=\frac{1}{\lambda_{\mu}}\,\Gamma^*_{\pm}(\mu,\mu), \nn \\
&\mL^{\pm}(\omega_{a})=\mT^{\pm}(\omega_{a})
\label{Lq01}
\end{align}
on any left-invariant 1-form $\omega_{a}$. 
The  contraction $\Gamma$ given in \eqref{Ggc} is said real and hermitian provided the conditions
\begin{align}
&\omega_{a}^*\wedge L(\omega_{b})=\Gamma(\omega_{a},\omega_{b})\mu,  \nn \\
&L(\omega^*)=(L(\omega_{a}))^*
\label{hrLq}
\end{align}
are satisfied. Such conditions amount to set
\beq
(\alpha\neq0,\,\beta=q^2\alpha,\,\nu,\,\gamma,\,\epsilon,\,\xi=\epsilon)\,\in\,\IR,
\label{hrc}
\eeq
together with  the implicit condition of invertibility $\Gamma_{\pm}(\mu,\mu)\neq 0$. 

The analysis of the spectra of the operators $\mL^{\pm}$ leads to characterise real and maximally hermitian contractions $\Gamma_{\pm}$ with respect to the Hodge operators $\mL^{\pm}$. 
Although long, straightforward calculations  show that the family of real and maximally hermitian 
contractions for $\mL^{\pm}$ \emph{does coincide} with that corresponding to $\mT^{\pm}$, whose parameters fullfill the conditions \eqref{as}. 

For real and maximally hermitian contractions it is further possible to  prove that the action of  $\mL^{\pm}(\xi)$ coincides with that of $\mT^{\pm}(\xi)$ on 2- and 3- left-invariant forms in \eqref{eiA}; this is indeed sufficient to prove that operators $\mL^{\pm}$ satisfy the same reality conditions satisfied by  $\mT^{\pm}$ (given in \eqref{rTt} and \eqref{rTm}),
\beq 
\lambda^{\pm}_{\xi^*}\,(\mL^{\pm}(\xi))^*=\lambda_{\xi}^{\pm}\mL^{\pm}(\xi^*).
\label{rcL}
\eeq
For this set of contractions, one can calculate that
\beq
\Gamma_{\pm}(\mu,\mu)=-2(q+q^{-1})^2(q^2+1+q^{-2})\,m^2\,\alpha\,\beta\,\epsilon^2.
\label{v4}
\eeq
Recalling \eqref{p15}, it is  then 
\beq
\mL^{\pm}(\mu)\neq\mT^{\pm}(\mu). 
\label{n4}
\eeq
This difference says that relations like \eqref{sTg}, \eqref{sTgm}, \eqref{fq} can not (for any choice of the real scale parameter $m$) be satisfied by the Hodge operators $\mL^{\pm}$. This is the reason why we do not extend them to proper Hodge dualities as in \eqref{shc}, but only use to give a different characterisation of the class of real and maximally hermitian  contractions.

\begin{rema}
A Hodge operator on $\Omega(\sq)$, the left covariant three dimensional exterior algebra on the homogeneous quantum Podle\'s sphere $\sq$ induced by the $\Omega(\SU)$ we have considered here, has been introduced in \cite{ag10} via an inner product obtained as a restriction to $\Omega(\sq)$ 
of a suitable class of contractions on $\Omega(\SU)$. 

It is immediate to see that the condition $\Gamma(\omega_{0},\omega_{0})=\nu=0$ \eqref{as} amounts to say  that the class of real and maximally hermitian contractions on $\Omega(\SU)$ would induce a degenerate inner product on $\Omega(\sq)$. The problem of defining a Hodge operator on $\Omega(\sq)$ whose corresponding Laplacian give a spectrum equivalent to the square of the Dirac operator \cite{sg10} requires further work.
\end{rema}


\appendix

\section{Exterior differential calculi  on Hopf algebras}\label{ass:a1}

 This appendix recalls, from the theory of differential calculi on algebras and  quantum groups,   definitions and general results needed to make the description of the  exterior algebras over  $\SU$ self consistent.  A more complete analysis is in  \cite{wor89,KS97}.

Let $\ch$ be the quantum group given by the unital $*$-Hopf algebra  $\ch=(\ch,\Delta,\varepsilon,S)$ over $\IC$, with  $\Omega^{1}(\ch)$ an $\ch$-bimodule. Given the  linear map $\dd:\ch\to\Omega^1(\ch)$, the pair $(\Omega^1(\ch), \dd)$ is a (first order) differential calculus over $\ch$  if $\dd$ satisfies the Leibniz rule, $\dd (h\, h^{\prime}) = (\dd  h) h^{\prime}  + h\, \dd h^{\prime} $ for $h,h^{\prime}\in \ch$, and 
if $\Omega^1(\ch)$ is generated by $\dd(\ch)$ as a $\ch$-bimodule. It is called  a $*$-calculus if there is an anti-linear involution $*:\Omega^1(\ch)\to\Omega^1(\ch)$ such that $(h_{1}(\dd h)h_{2})^*=h_{2}^{*}(\dd(h^*))h_{1}^{*}$ for any $h,h_{1},h_{2}\in \ch$. 

A first order differential calculus is said left covariant provided  a left coaction $\Delta_{L}^{(1)}:\Omega^1(\ch)\to\ch\otimes\Omega^1(\ch)$ exists, such that $\Delta_{L}^{(1)}(\dd h)=(1\otimes \dd)\Delta(h)$ and $\Delta_{L}^{(1)}(h_1\,\alpha\,h_2)=\Delta(h_{1})\Delta_{L}^{(1)}(\alpha)\Delta(h_2)$ for any $h,h_1,h_2\in\,\ch$ and $\alpha\in\,\Omega^1(\ch)$.

The property of right covariance is stated in complete analogy; 
there is in addition the notion of a bicovariant calculus, namely a calculus which is both left and right covariant 
and satisfying the  compatibility condition:
$$
(1\otimes\Delta^{(1)}_{R})\circ\Delta^{(1)}_{L}=(\Delta^{(1)}_{L}\otimes 1)\circ\Delta^{(1)}_{R}.
$$
Left covariant calculi on $\ch$ are  in one to one correspondence with  right ideals 
$\mathcal{Q} \subset \ker\eps$. There is a left $\ch$-modules isomorphism given by 
$\Omega^{1}(\ch)\simeq\ch\otimes(\ker\eps/\cq)$, which allows to  recover the complex vector space $\ker\eps/\cq$ as the set of left invariant 1-forms, namely the elements $\omega_{a}\in \Omega^{1}(\ch)$ such that 
$$
\Delta^{(1)}_{L}(\omega_{a})=1\otimes\omega_{a} .
$$
The dimension of $\ker\eps/\cq$ is referred to as the dimension of the calculus. 
A left covariant first order differential calculus is  a $*$-calculus if and only if $(S(Q))^*\in \cq$ for any $Q\in \cq$. If this is the case,  the left coaction of $\ch$ on $\Omega^{1}(\ch)$ is compatible with the $*$-structure: 
$\Delta^{(1)}_{L}(\dd h^{*})=(\Delta^{(1)}(\dd h))^*$.
Bicovariant calculi  corresponds to   right ideals $\mathcal{Q}\subset\ker\varepsilon$ which are in addition stable under the right adjoint coaction $\Ad$ of $\ch$ onto itself, that is to say $\Ad(\cq) \subset \cq \otimes \ch$\footnote{  
Explicitly, $\Ad=\left(\id \otimes m \right) \left(\tau\otimes
\id \right)\left(S\otimes\Delta\right)\Delta$, with $\tau$ the flip
operator, or $\Ad(h) = \co{h}{2} \otimes \left(S(\co{h}{1}) \co{h}{3} \right)$
in Sweedler notation.}. With $G$ a compact connected and real Lie group, and $\ch=\ca(G)$ the standard Hopf $*$-algebra of complex valued functions defined on $G$, the well known natural bicovariant calculus  corresponds to the choice $\cq=(\ker\,\varepsilon)^2$. 

The right coaction of $\ch$ on  $\Omega^{1}(\ch)$ defines matrix elements 
\beq \label{ri-co-form}
\Delta^{(1)}_{R}(\omega_{a})=\sum_b \omega_{b}\otimes J_{ba}.
\eeq
where $J_{ab}\in\ch$. This matrix is invertible, since $\sum_b S(J_{ab})J_{bc}=\delta_{ac}$ and 
$\sum_b J_{ab}S(J_{bc})=\delta_{ac}$. In addition one finds that $\Delta(J_{ab})=\sum_c J_{ac}\otimes J_{cb}$ and 
$\varepsilon(J_{ab})=\delta_{ab}$.  It gives a basis of right invariant 1-forms, $\eta_{a}=\omega_{b}S(J_{ba})$

In order to construct an exterior algebra $\Omega(\ch)$ over the bicovariant first order differential calculus $(\Omega^{1}(\ch),\dd)$ one uses a  braiding map replacing the flip automorphism $\tau$. Define the tensor product $\Omega^{\otimes k}(\ch)=\Omega^{1}(\ch)\otimes_{\ch}\ldots\otimes_{\ch}\Omega^{1}(\ch)$
with $k$ factors. 
There exists a unique $\ch$-bimodule homomorphism 
$\sigma:\Omega^{\otimes 2}(\ch)\to\Omega^{\otimes 2}(\ch)$ 
such that $\sigma(\omega\otimes\eta)=\eta\otimes\omega$ for any left invariant 1-form $\omega$ and any right invariant 1-form $\eta$. The map $\sigma$ is invertible and commutes with the left coaction of $\ch$: 
$$
(1\otimes\sigma)\circ\Delta_{L}^{(2)}=\Delta_{L}^{(2)}\circ\sigma,
$$ 
with $\Delta_{L}^{(2)}$ the extension of the coaction to the tensor product. There is an analogous invariance for the right 
coaction. 
Moreover, $\sigma$ satisfies a braid equation. On $\Omega^{\otimes 3}(\ch)$: 
\beq
(1\otimes\sigma)\circ(\sigma\otimes 1)\circ(1\otimes \sigma)=(\sigma\otimes 1)\circ(1\otimes \sigma)\circ(\sigma\otimes 1).
\label{br3}
\eeq
The braiding map provides  an antisymmetrizer operator
$A^{(k)}:\Omega^{\otimes k}(\ch)\to\Omega^{\otimes k}(\ch)$ for any $k\geq2$. 
Consider the permutation group $P(k)$ of $k$ elements. Any $p\in P(k)$ can be written as a product $p=\Pi_{j\in I(p)}\tau_{j_{1}}\cdots\tau_{j_{\natural( I(p))}}$ of $I(p)$ ($\natural(I(p))$ is the cardinality of $I(p)$) nearest neighbour transpositions $\{\tau_{1},\ldots,\tau_{k-1}\}$ (i.e. elementary flips), where  $\tau_{j}$ transposes $j$ with $j+1$ leaving all other elements in $\{1,\ldots,k\}$ fixed. Using the braiding, set $\sigma_{j}:\Omega^{\otimes k}(\ch)\to\Omega^{\otimes k}(\ch)$ as $\sigma_{j}:=1\otimes\ldots\otimes\sigma\otimes\ldots\otimes 1$, where the product contains $k-1$ factors and $\sigma$ occurs in the $j^{th}$ place. It is then immediate to associate  $P(k)\ni\,p\,\mapsto\,\mathfrak{P}=\Pi_{j\in I(p)}\sigma_{j_{1}}\cdots\sigma_{j_{\natural I(p)}}$, and  to define 
\beq
\label{Ank}
A^{(k)}:=\Pi_{p\in P(k)}\,\mathrm{sign}(p)\,\mathfrak{P}(p)
\eeq
as the (order $k$) antisymmetriser operator associated to the braiding of the calculus. It is a bimodule automorphism; 
the Hopf ideals $\mathcal{S}_{\cq}^{(k)}=\ker\,A^{(k)}$ give the quotients
\beq
\Omega^{k}(\ch)=\Omega^{1}(\ch)^{\otimes k}/\mathcal{S}^{(k)}_{\cq}
\label{wedk}
\eeq
the structure of a $\ch$-bicovariant bimodule  which can be written as $\Omega^{k}(\ch)=\mathrm{Range}\,A^{(k)}$. The exterior algebra is  
$(\Omega(\ch)=\oplus_{k}\Omega^{k}(\ch),\wedge)$
with the identification $\Omega^{0}(\ch)=\ch$.  The exterior derivative is extended to $\Omega(\ch)$ as the only degree one derivation such that $\dd^2=0$. The  algebra $\Omega(\ch)$ has natural left and right $\ch$-comodule structure, given by recursively setting
$$
\Delta_{L}^{(k+1)}(\dd\theta)=(1\otimes\dd)\Delta_{L}^{(k)}(\theta), \qquad \Delta_{R}^{(k+1)}(\dd\theta)=(\dd\otimes 1)\Delta_{R}^{(k)}(\theta).
$$
Finally, the $*$-structure on $\Omega^{1}(\ch)$ is  extended to an antilinear $*:\Omega(\ch)\to\Omega(\ch)$ by $(\theta\wedge
 \theta^{\prime})^{*}=(-1)^{kk^{\prime}}\theta^{\prime*}\wedge\theta^{*}$ with $\theta\in \Omega^{k}(\ch)$ and $\theta^{\prime}\in \Omega^{k^{\prime}}(\ch)$; the exterior derivative operator satisfies the identity $(\dd\theta)^{*}=\dd(\theta^*)$. 

The square of the braiding map does not in general coincide  with the identity.  The inverse map $\sigma^{-1}$ also satisfies a braid equation as in \eqref{br3}, and provides, in analogy to \eqref{Ank}, a set of antisymmetrizer operators   $A^{(k)}_{-}: \Omega^{1}(\ch)^{\otimes k}\to\Omega^{1}(\ch)^{\otimes k}$, giving Hopf ideals 
$\mathcal{S}_{\mathcal{Q}-}^{(k)}=\ker\,A^{(k)}=\ker\,A^{(k)}=\mathcal{S}_{\cq}$.
The quotients $\Omega^k_-(\ch)=\Omega^1(\ch)^{\otimes k}/\mathcal{S}_{\cq-}^{(k)}$ describe then an exterior algebra 
$(\Omega_-(\ch)=\oplus_{k}\Omega_{-}^{k}(\ch),\vee)$ which is isomorphic to the $\Omega(\ch)$ 
considered above.

\subsection{Derivations associated to exterior differential calculi}

The tangent space of
the calculus is  the complex vector space of elements  out of $\ch^{\prime}$ -- the dual space $\ch^{\prime}$ of functionals on $\ch$ --  defined by 
$$
\mathcal{X}_{\mathcal{Q}}:=
\{L\in \ch^{\prime} ~:~ L(1)=0,\,L(Q)=0, \,\, \forall \, Q\in\mathcal{Q}\}. 
$$
There exists a unique bilinear form 
\beq
\hp{~}{~}:\mathcal{X}_{\cq}\times\Omega^{1}(\ch), \qquad\hp{L}{x \dd y}:=\eps(x)L(y) ,
\label{dptg}
\eeq
giving a non-degenerate dual pairing between the vector spaces $\mathcal{X}_{\cq}$ and $\ker\eps/\cq$, and then a vector space isomorphism $\mathcal{X}_{\cq}\simeq(\ker\varepsilon/\cq)$. 
The dual space $\ch^{\prime}$ has natural left and right (mutually commuting) actions  on $\ch$:  
\beq
L\triangleright h:=h_{(1)}X(h_{(2)}),\qquad
 h\triangleleft L:=X(h_{(1)})h_{(2)}.
\label{deflr}
\eeq
If the  vector space  $\mathcal{X}_{\cq}$ is finite dimensional  its elements belong to the dual Hopf algebra 
$ \ch^{\prime} \supset \ch^{o} = (\ch^{o}, \Delta_{\ch^o}, \eps_{\ch^o},S_{\ch^o})$, defined as the largest Hopf $*$-subalgebra contained in $\ch^{\prime}$. 
In such a case  the 
$*$-structures are compatible with both actions,
$$
L \lt h^* = ((S(L))^* \lt h)^*,\qquad
h^* \rt  L = (h \rt  (S(L))^*)^*,
$$
for any $ L \in \ch^{o}, \ h \in \ch$ 
and the exterior derivative can be written as:
\beq
\dd h := \sum_a
~(L_{a} \triangleright h) ~\omega_{a} =\sum_{a} \omega_{a}(-S^{-1}(L_{a}))\lt h, 
\label{ded}
\eeq
where $\hp{L_{a}}{\omega_{b}}=\delta_{ab}$.   
The twisted Leibniz rule of derivations of the basis elements $L_{a}$ is dictated by their coproduct:  
\beq
\Delta_{\ch^o}(L_{a})=1\otimes L_{a}+
\sum\nolimits_b L_{b}\otimes f_{ba},
\label{cpuh}
\eeq
where the $f_{ab} \in \ch^{o}$ consitute an algebra representation of $\ch$, also  controlling the $\ch$-bimodule structure of $\Omega^{1}(\ch)$:
\begin{align}\label{bi-struct}
\omega_{a} h = \sum\nolimits_b (f_{ab}\triangleright h)\omega_{b} , \qquad 
h \omega_{a} = \sum\nolimits_b \omega_{b} \left( (S^{-1}(f_{ab}) )\triangleright h \right) , 
\qquad \mathrm{for} \quad h \in \ch. 
\end{align}

\section{The quantum group $\SU$}
\label{eSU}
The polynomial algebra $\ASU$ of the quantum group
$\SU$  is the unital  $*$-algebra generated by elements $a$ and $c$, with relations
\begin{align}
\label{derel}
& ac=qca\quad ac^*=qc^*a\quad cc^*=c^*c , \nn \\
& a^*a+c^*c=aa^*+q^{2}cc^*=1 .
\end{align}
The deformation parameter
$q\in\IR$ can be restricted without loss of generality to the interval $0<q<1$:  in the formal limit $q \to1$ one recovers the commutative coordinate algebra on the group manifold $\mathrm{SU(2)}$.
Relations \eqref{derel} are equivalent to the unitarity of the matrix  
$$ U = 
\left(
\begin{array}{cc} a & -qc^* \\ c & a^*
\end{array}\right) , 
$$
while  the Hopf algebra structures for $\ASU$ can be expressed as:
$$
\Delta\, U = U \otimes U ,  \qquad S(U) = U^* , \qquad \eps(U) = 1. 
$$
The quantum universal envelopping algebra $\su$ is the unital Hopf $*$-algebra
generated as an algebra by four elements $K^{\pm 1},E,F$, with $K K^{-1}=1$ and the
relations: 
\beq 
K^{\pm}E=q^{\pm}EK^{\pm}, \qquad 
K^{\pm}F=q^{\mp}FK^{\pm}, \qquad  
[E,F] =\frac{K^{2}-K^{-2}}{q-q^{-1}} . 
\label{relsu}
\eeq 
The $*$-structure is
$K^*=K, \,  E^*=F $,
and the Hopf algebra structure is provided  by coproduct
$$\Delta(K^{\pm}) =K^{\pm}\otimes K^{\pm}, \quad
\Delta(E) =E\otimes K+K^{-1}\otimes E,  \quad 
\Delta(F)
=F\otimes K+K^{-1}\otimes F,$$
while the antipode and counit are
\begin{align*}
&S(K) =K^{-1}, \quad
S(E) =-qE, \quad 
S(F) =-q^{-1}F \\ 
&\varepsilon(K)=1, \quad \varepsilon(E)=\varepsilon(F)=0.
\end{align*}
 The quadratic element 
\beq
C_{q}=\frac{qK^2-2+q^{-1}K^{-2}}{(q-q^{-1})^2}+FE-\tfrac{1}{4}
\label{cas}
\eeq
is a quantum Casimir operator that generates the centre of $\su$.

The Hopf $*$-algebras  $\su$ and $\ASU$ are dually paired. The $*$-compatible bilinear mapping $\hs{~}{~}:\su\times\ASU\to\IC$ is  given on the generators by
\begin{align}
&\langle K^{\pm},a\rangle=q^{\mp 1/2}, \qquad  \langle K^{\pm},a^*\rangle=q^{\pm 1/2},   \nn\\
&\langle E,c\rangle=1, \qquad \langle F,c^*\rangle=-q^{-1}, \label{ndp}
\end{align}
with all other couples of generators pairing to zero. This pairing is proved \cite{KS97} to be non-degenerate. The algebra $\su$ is recovered as a $*$-Hopf subalgebra in the dual algebra $\ASU^o$, the largest Hopf $*$-subalgebra contained in the dual vector space  $\ASU^{\prime}$.
There are \cite{wor87} $*$-compatible canonical commuting actions of $\su$ on $\ASU$: 
$$
h \lt x := \co{x}{1} \,\hs{h}{\co{x}{2}}, \qquad
x \rt  h := \hs{h}{\co{x}{1}}\, \co{x}{2}. 
$$
Given the algebra $\ca(\U(1)):=\IC[z,z^*] \big/ \!\!<zz^* -1>$,  the map  
\beq  \label{qprp}
\pi: \ASU \, \to\,\ca(\U(1)) , \qquad  
\pi\,\left(
\begin{array}{cc} a & -qc^* \\ c & a^*
\end{array}\right):=
\left(
\begin{array}{cc} z & 0 \\ 0 & z^*
\end{array}\right)
\eeq 
is a surjective Hopf $*$-algebra homomorphism, so that  $\U(1)$
is a quantum subgroup of $\SU$ with right coaction:
\beq 
\delta_{R}:= (\id\otimes\pi) \circ \Delta \, : \, \ASU \,\to\,\ASU \otimes
\ca(\U(1)) . \label{cancoa} 
\eeq 
This right coaction gives a decomposition 
\beq
\label{clnm}
\ASU=\oplus_{n\in\,\IZ}\cl_{n}, \qquad\qquad\cl_{n}:=\{x\in\,\ASU\;:\;\delta_{R}(x)=x\otimes z^{-n}\},
\eeq
with $\Asq=\cl_{0}$ the algebra of the standard Podle\'s sphere, each $\Asq$-bimodule $\cl_{n}$ gives the set of (charge $n$)  $\U(1)$-coequivariant maps for the topological quantum principal bundle $\Asq\hookrightarrow\ASU$.

\subsection{A 4D exterior algebra over the quantum group $\SU$}\label{se:4dc} 

This section presents, following \cite{ag10},  the exterior algebra over the  so called $4D_{+}$ bicovariant calculus on $\SU$, which was introduced  as a first order differential calculus  in \cite{wor89}. 

The ideal $\cq_{\SU}\subset\ker\varepsilon_{\SU}$ corresponding to the $4D_+$ calculus is generated by the nine elements $\{c^{2}; \,  c(a^{*}-a); \, q^{2}a^{*2}-(1+q^{2})(aa^{*}-cc^{*})+a^{2}; \, c^{*}(a^{*}-a); \, c^{*2}; \, [q^{2}a+a^{*}-q^{-1}(1+q^{4})]c; \, 
[q^{2}a+a^{*}-q^{-1}(1+q^{4})](a^{*}-a); \, [q^{2}a+a^{*}-q^{-1}(1+q^{4})]c^{*}; \, [q^{2}a+a^{*}-q^{-1}(1+q^{4})][q^{2}a+a^{*}-(1+q^{2})]\}$. One has $\mathrm{Ad}(\cq_{\SU})\subset\cq_{\SU}\otimes\ASU$ and $\dim(\ker\varepsilon_{\SU}/{\cq_{\SU}})=4$. The quantum tangent space is a  four dimensional $\mathcal{X}_{\cq}\subset\su$, with a basis given by
\begin{align}
& L_{-}=q^{\frac{1}{2}}FK^{-1},\qquad L_{z}=\frac{K^{-2}-1}{q-q^{-1}},\qquad L_{+}=q^{-\frac{1}{2}}EK^{-1}; \nn \\
& L_{0}=\frac{q(K^{2}-1)+q^{-1}(K^{-2}-1)}{(q-q^{-1})^{2}}\,+FE=
\frac{q(K^{-2}-1)+q^{-1}(K^{2}-1)}{(q-q^{-1})^{2}}\,+EF.
\label{Lq}
\end{align}
The vector $L_{0}$ belongs to the centre of $\su$, and  differs from the quantum Casimir \eqref{cas} by a constant term,
\beq\label{casbis}
C_{q}=L_{0}+[\tfrac{1}{2}]^{2}-\tfrac{1}{4}.
\eeq
The differential $\dd:\ASU\mapsto\Omega^{1}(\SU)$ is written for any $x\in \ASU$ as (with $R_{a}=-S^{-1}(L_{a})$)
\beq
\dd x=\sum_a (L_{a}\lt x)\omega_{a}=\sum_a \omega_{a}(R_{a}\lt x), \qquad\qquad \mathrm{with}\,\, R_{a}=-S^{-1}(L_{a})
\label{d4}
\eeq
on the dual  basis of left invariant forms $\omega_{a}\in\Omega^{1}(\SU)$,  which are  
\begin{align}
&\omega_{-}=c^{*}\dd a^{*}-qa^{*}\dd c^{*}, \nn \\
&\omega_{+}=a\dd c-qc\dd a, \nn \\ 
&\omega_{z}=a^{*}\dd a+c^{*}\dd c-(a\dd a^{*}+q^{2}c\dd c^{*}),\nn \\
&\omega_{0}=(1+q)^{-1}\rho^{-1}[a^{*}\dd a+c^{*}\dd c+q(a\dd a^{*}+q^{2}c\dd c^{*})],
\label{om4}
\end{align}
where $\rho=[\frac{1}{2}][\frac{3}{2}]$. This differential calculus reduces in the classical limit  to the standard  three-dimensional bicovariant calculus on $\mathrm{SU(2)}$, since one has $\omega_0\to0$ for $q\to1$. Since $(S(\cq_{\SU}))^*\subset\cq_{\SU}$,  the $*$-structure on $\ASU$ is  consistently extended to an antilinear  $*$-structure on $\Omega^{1}(\SU)$, such that $(\dd x)^*=\dd (x^*)$ for any $x\in \ASU$. For the basis of left invariant 1-forms it is
\beq
\omega_{-}^{*}=-\omega_{+},\qquad\omega_{z}^{*}=-\omega_{z},\qquad\omega_{0}^{*}=-\omega_{0}.
\label{ss}
\eeq  

The construction of the exterior algebras $\Omega(\SU)\simeq\Omega_{-}(\SU)$ is given via the braiding map $\sigma:\Omega^{\otimes2}(\SU)\to\Omega^{\otimes2}(\SU)$ and its inverse. 
\begin{align}
\sigma^{\pm}&(\omega_{-}\otimes\omega_{-})=\omega_{-}\otimes\omega_{-},  \label{sigmm} \\ \nn \\
\sigma&\left(\begin{array}{c} \omega_-\otimes\omega_{0} \\ \omega_0\otimes\omega_{-} \\ \omega_-\otimes\omega_{z} \\ \omega_{z}\otimes\omega_{-} \end{array}\right)=\left(\begin{array}{cccc}
1-q^2 & 1 & 0 & 0 \\ q^2 & 0 & 0 & 0 \\  1+q^2 & 0 & 0 & 1 \\ -1-q^{-2} & 0 & q^{-2} & 1-q^{-2} \end{array}
\right) \left(\begin{array}{c} \omega_-\otimes\omega_{0} \\ \omega_0\otimes\omega_{-} \\ \omega_-\otimes\omega_{z} \\ \omega_{z}\otimes\omega_{-} \end{array}\right), \nn \\
\sigma^{-1}&\left(\begin{array}{c} \omega_-\otimes\omega_{0} \\ \omega_0\otimes\omega_{-} \\ \omega_-\otimes\omega_{z} \\ \omega_{z}\otimes\omega_{-} \end{array}\right)=\left(\begin{array}{cccc}
0 & q^{-2} & 0 & 0 \\ 1 & 1-q^{-2} & 0 & 0 \\  0 & 1+q^2 &  1-q^2 & q^2 \\ 0 & -1-q^{-2} & 1 & 0 \end{array}
\right) \left(\begin{array}{c} \omega_-\otimes\omega_{0} \\ \omega_0\otimes\omega_{-} \\ \omega_-\otimes\omega_{z} \\ \omega_{z}\otimes\omega_{-} \end{array}\right),  \label{sigmz} 
\end{align}
\begin{align}
\sigma&\left(\begin{array}{c} \omega_z\otimes\omega_{0} \\ \omega_0\otimes\omega_{z} \\ \omega_z\otimes\omega_{z} \\ \omega_{-}\otimes\omega_{+} \\ \omega_{+}\otimes\omega_{-} \end{array}\right)=\left(\begin{array}{ccccc}
-(q-q^{-1})^2 & 1 & 0 & -(q-q^{-1})^2 & (q-q^{-1})^2 \\ 1 & 0 & 0 & 0 & 0 \\   q^2-q^{-2} & 0 &   1 & q^2-q^{-2}  & -(q^2-q^{-2})  \\ -1  & 0 & 0 & 0 & 1 \\ 1 & 0 & 0 & 1 & 0 \end{array}
\right) \left(\begin{array}{c} \omega_z\otimes\omega_{0} \\ \omega_0\otimes\omega_{z} \\ \omega_z\otimes\omega_{z} \\ \omega_{-}\otimes\omega_{+} \\ \omega_{+}\otimes\omega_{-}  \end{array}\right), \nn \\
\sigma^{-1}&\left(\begin{array}{c} \omega_z\otimes\omega_{0} \\ \omega_0\otimes\omega_{z} \\ \omega_z\otimes\omega_{z} \\ \omega_{-}\otimes\omega_{+} \\ \omega_{+}\otimes\omega_{-} \end{array}\right)=\left(\begin{array}{ccccc}
0 & 1 & 0 & 0 & 0 \\ 1 & - (q-q^{-1})^2 & 0  & -(q-q^{-1})^2 & (q-q^{-1})^2 \\  0 &  q^2-q^{-2} & 1 &   q^2-q^{-2}  & -(q^2-q^{-2})  \\ 0  & -1 & 0 & 0 & 1 \\ 0 & 1 & 0 & 1 & 0 \end{array}
\right) \left(\begin{array}{c} \omega_z\otimes\omega_{0} \\ \omega_0\otimes\omega_{z} \\ \omega_z\otimes\omega_{z} \\ \omega_{-}\otimes\omega_{+} \\ \omega_{+}\otimes\omega_{-}  \end{array}\right),  \label{sigzz} \\ \nn \\
\sigma^{\pm}&(\omega_{0}\otimes\omega_{0})=\omega_{0}\otimes\omega_{0}  \label{sig00} 
\end{align}
\begin{align}
\sigma&\left(\begin{array}{c} \omega_+\otimes\omega_{0} \\ \omega_0\otimes\omega_{+} \\ \omega_+\otimes\omega_{z} \\ \omega_{z}\otimes\omega_{+} \end{array}\right)=\left(\begin{array}{cccc}
1-q^{-2} & 1 & 0 & 0 \\ q^{-2} & 0 & 0 & 0 \\  -1-q^{-2} & 0 & 0 & 1 \\ 1+q^{2} & 0 & q^{2} & 1-q^{2} \end{array}
\right) \left(\begin{array}{c} \omega_+\otimes\omega_{0} \\ \omega_0\otimes\omega_{+} \\ \omega_+\otimes\omega_{z} \\ \omega_{z}\otimes\omega_{+} \end{array}\right), \nn \\ 
\sigma^{-1}&\left(\begin{array}{c} \omega_+\otimes\omega_{0} \\ \omega_0\otimes\omega_{+} \\ \omega_+\otimes\omega_{z} \\ \omega_{z}\otimes\omega_{+} \end{array}\right)=\left(\begin{array}{cccc}
0 & q^2 & 0 & 0 \\ 1 & 1-q^{2}  & 0 & 0 \\  0 & -1-q^{-2} & 1-q^{-2}  & q^{-2} \\ 0 & 1+q^{2} & 1 & 0  \end{array}
\right) \left(\begin{array}{c} \omega_+\otimes\omega_{0} \\ \omega_0\otimes\omega_{+} \\ \omega_+\otimes\omega_{z} \\ \omega_{z}\otimes\omega_{+} \end{array}\right),  \label{sigpz} \\ \nn \\ 
\sigma&(\omega_{+}\otimes\omega_{+})=\omega_{+}\otimes\omega_{+}. \label{sigpp}
\end{align}
The spectral resolution of the braiding gives
\beq
(\sigma^{\pm}+1)(\sigma^{\pm}+q^{2})(\sigma^{\pm}+q^{-2})=0;
\label{spb}
\eeq
it allows to give the spectral resolution of the antisymmetriser operators \cite{sch99}. With $\Omega^{k}(\SU)\,\ni\,\omega_{a_{1}}\wedge\ldots\wedge\omega_{a_{k}}=A^{(k)}(\omega_{a_{1}}\otimes\ldots\otimes\omega_{a_{k}})$, one has from \cite{ag10} that: 
\begin{align}
\varphi_{\pm}:=\omega_{\mp}\wedge\omega_{0},\qquad&\qquad A^{(2)}(\varphi_{\pm})=(1+q^{\pm2})\varphi_{\pm} \nn \\
\kappa_{\pm}:=\omega_{\pm}\wedge\omega_{0}\mp(1-q^{\pm2})\omega_{\pm}\wedge\omega_{z}, \qquad&\qquad A^{(2)}(\kappa_{\pm})=(1+q^{\pm2})\kappa_{\pm} \nn \\
\psi_{\pm}:=\omega_{0}\wedge\omega_{z}+(1-q^{\pm 2})\omega_{-}\wedge\omega_{+},\qquad&\qquad A^{(2)}(\psi_{\pm})=(1+q^{\pm2})\psi_{\pm}
\label{spe2}
\end{align}
on a basis of  the 6-dimensional vector space of left-invariant 2-forms, with $\varphi^*_{\pm}=\varphi_{\mp}, \,\kappa^*_{\pm}=\kappa_{\mp},\,\psi_{\pm}^*=\psi_{\mp}$. On 3-forms it is
\beq
\label{p13}
\begin{array}{lll}
\chi_{-}:=\omega_{+}\wedge\omega_{0}\wedge\omega_{z},&  &
\chi_{+}:=\omega_{-}\wedge\omega_{0}\wedge\omega_{z} \\
\chi_{0}:=\omega_{-}\wedge\omega_{+}\wedge\omega_{z}, & & 
\chi_{z}:=\omega_{-}\wedge\omega_{+}\wedge\omega_{0},
\end{array}
\eeq
\beq
A^{(3)}(\chi_{a})=2(1+q^2+q^{-2})\chi_{a} 
\label{p15-}
\eeq
for $a=-,+,z,0$, with $\chi_{-}^*=-q^{-2}\chi_{+}$, $\chi_{0}^*=\chi_{0}$ and $\chi_{z}^*=\chi_{z}$. The exterior algebra corresponding to this calculus is 4-dimensional; the $\ASU$-bimodule $\Omega^{4}(\SU)$  is 1 dimensional. The action of the antisymmetriser on the left invariant basis element $\vartheta=\omega_{-}\wedge\omega_{+}\wedge\omega_{z}\wedge\omega_{0}$ gives 
\beq
\label{p15}
A^{(4)}(\vartheta)=2(q^4+2q^2+6+2q^{-2}+q^{-4})\vartheta.
\eeq
The analysis of the exterior algebra $\Omega_{-}(\SU)$, introduced in appendix \S \ref{ass:a1} and constructed via the "inverse braiding" $\sigma^{-1}$ map,  shows a similar pattern. With $\Omega_{-}^{k}(\SU)\,\ni\,\omega_{a_{1}}\vee\ldots\vee\omega_{a_{k}}=A_{-}^{(k)}(\omega_{a_{1}}\otimes\ldots\otimes\omega_{a_{k}})$, one has
\begin{align}
\check{\varphi}_{\pm}:=\omega_{\mp}\vee\omega_{0}=q^{\mp2}\varphi_{\pm},\qquad&\qquad A_{-}^{(2)}(\check{\varphi}_{\pm})=(1+q^{\mp2})\check{\varphi}_{\pm}
\nn \\ \check{\kappa}_{\pm}:=\omega_{\pm}\vee\omega_{0}\mp(1-q^{\pm2})\omega_{\pm}\vee\omega_{z}=q^{\mp2}\kappa_{\pm},\qquad&\qquad A_{-}^{(2)}(\check{\kappa}_{\pm})=(1+q^{\mp2})\check{\kappa}_{\pm} \nn \\
\check{\psi}_{\pm}:=\omega_{0}\vee\omega_{z}+(1-q^{\pm2})\omega_{-}\vee\omega_{+}=q^{\mp2}\psi_{\pm},\qquad&\qquad A_{-}^{(2)}(\check{\psi}_{\pm})=(1+q^{\mp2})\check{\psi}_{\pm}
\label{spe2-}
\end{align}
on 2-forms, with $\check{\varphi}^*_{\pm}=\check{\varphi}_{\mp}, \,\check{\kappa}^*_{\pm}=\check{\kappa}_{\mp},\,\check{\psi}_{\pm}^*=\check{\psi}_{\mp}$;  while on 3- and 4-forms, with the same notations used in \eqref{p13},
\begin{align}
\check{\chi}_{a}=\chi_{a},\qquad&\qquad A_{-}^{(3)}(\check{\chi}_{a})=2(1+q^{2}+q^{-2})\check{\chi}_{a} \nn \\
\check{\vartheta}=\omega_{-}\vee\omega_{+}\vee\omega_{0}\vee\omega_{z}=\vartheta, \qquad&\qquad
A_{-}^{(4)}(\check{\vartheta})=2(q^4+2q^2+6+2q^{-2}+q^{-4})\check{\vartheta}.
\label{sp3-}
\end{align}
These relations show that the anti-linear $*$-structure  does not commute 
with the action of the antisymmetrisers $A_{\pm}^{(2)}$.

\bigskip

\subsection*{Acknowledgements} I should like to thank Giovanni Landi: he gave me the chance to  evolve in a common research project  those  themes that I found more interesting, indeed giving me a precious guidance at each stage of the work. I should like to thank Lucio Cirio for his suggestions, and Istvan Heckenberger for his constant feedback and comments. 

This paper closes a long research period I spent in Paris and in Bonn. It is a pleasure to acknowledge the support of the I.H.E.S. (Bures sur Yvette), the Max-Planck-Institut f\"ur Mathematik  in Bonn, the Hausdorff Zentrum f\"ur Mathematik in Bonn. And I should like to warmly thank Sergio Albeverio, Yuri I. Manin, Matilde Marcolli: they were always at my side during the  years in Rheinland.
 
This paper has been then revised while working at the Mathematisches Institut der L.M.U., M\"unchen. I should like to express my gratitude to Detlef D\"urr for his mentorship and  support.

\end{document}